\documentclass[11pt]{article}
\usepackage{amsmath,amssymb,amsthm,morefloats}
\usepackage{epsfig,psfrag,url}
\usepackage{graphicx}
\usepackage{verbatim}
\usepackage{subfigure}

\numberwithin{equation}{section}

\newcommand{\goto}{\rightarrow}
\newcommand{\bigo}{{\mathcal O}}

\def\Xint#1{\mathchoice
   {\XXint\displaystyle\textstyle{#1}}%
   {\XXint\textstyle\scriptstyle{#1}}%
   {\XXint\scriptstyle\scriptscriptstyle{#1}}%
   {\XXint\scriptscriptstyle\scriptscriptstyle{#1}}%
   \!\int}
\def\XXint#1#2#3{{\setbox0=\hbox{$#1{#2#3}{\int}$}
     \vcenter{\hbox{$#2#3$}}\kern-.5\wd0}}

\def\dashint{\Xint-}

\DeclareMathOperator{\imag}{Im}

\DeclareMathOperator{\sign}{sign}

\newenvironment{choices}{\left\{ \begin{array}{ll}}{\end{array}\right.}
\newcommand\when{&\text{if~}}
\newcommand\otherwise{&\text{otherwise}}

\newenvironment{mat}{\left[\begin{array}{ccccccccccccccc}}{\end{array}\right]}
\newcommand\bcm{\begin{mat}}
\newcommand\ecm{\end{mat}}

\newcommand{\bea}{\begin{eqnarray}}
\newcommand{\eea}{\end{eqnarray}}
\newcommand{\bean}{\begin{eqnarray*}}
\newcommand{\eean}{\end{eqnarray*}}
\newcommand{\ba}{\begin{array}}
\newcommand{\ea}{\end{array}}
\newcommand{\beqs}{\begin{equation*}\begin{split}}

\newtheorem{remark}{Remark}[section]

\newtheorem{lemma}{Lemma}[section]

\newtheorem{theorem}{Theorem}[section]

\newtheorem{definition}{Definition}[section]

\newtheorem{proposition}{Proposition}[section]

\long\def\symbolfootnote[#1]#2{\begingroup%
\def\thefootnote{\fnsymbol{footnote}}\footnote[#1]{#2}\endgroup}

\newcommand{\D}{\text{d}}

\begin{document}
\title{Rational approximation, oscillatory Cauchy integrals and Fourier transforms}
\author{Thomas Trogdon$^1$\\
\phantom{.}\\
Courant Institute of Mathematical Sciences\\
New York University\\
251 Mercer St\\
New York, NY 10012, USA\\}
\maketitle

\footnotetext[1]{Email: trogdon@cims.nyu.edu}

\begin{abstract}
We develop the convergence theory for a well-known method for the interpolation of functions on the real axis with rational functions. Precise new error estimates for the interpolant are derived using existing theory for trigonometric interpolants.  Estimates on the Dirichlet kernel are used to derive new bounds on the associated interpolation projection operator.  Error estimates are desired partially due to a recent formula of the author for the Cauchy integral of a specific class of so-called oscillatory rational functions.  Thus, error bounds for the approximation of the Fourier transform and Cauchy integral of oscillatory smooth functions are determined.  Finally, the behavior of the differentiation operator is discussed.  The analysis here can be seen as an extension of that of Weber (1980) and Weideman (1995) in a modified basis used by Olver (2009) that behaves well with respect to function multiplication and differentiation.
\end{abstract}

\section{Introduction}

Trigonometric interpolation with the discrete Fourier transform is a classic approximation theory topic and there exists a wide variety of results.  See \cite{atkinson} for an in-depth discussion.  The ``fast'' nature of fast Fourier transform (FFT) makes this type of approximation appealing.  The FFT is just as easily considered as a method to compute a Laurent expansion of function on a circle in the complex plane centered at the origin \cite{Henrici1979}.  Furthermore, once a function is expressed in a Laurent expansion (assuming sufficient decay of the Laurent coefficients) it may be mapped to a function on the real axis with a M\"obius transformation.  This idea has been exploited many times.  It was used in \cite{Henrici1979} to compute Laplace transforms, in \cite{Weber1980} to compute Fourier transforms, in \cite{SOHilbertTransform,Weideman1995} to compute Hilbert transforms, in \cite{Trogdon2013} to compute oscillatory singular integrals and in \cite{SOPainleveII,Trogdon2013} to solve Riemann--Hilbert problems.  Specifically, the method discussed here can be realized by using interpolation on $[0,2\pi)$ with complex exponentials and composing the interpolant with a map from real axis to the unit circle (an $\arctan$ transformation).  This produces a rational approximation of a function on $\mathbb R$.

The goal of the current paper is three-fold.   First, we follow the convergence theory of interpolants on the periodic interval $\mathbb T = [0,2\pi)$ (see \cite{Kress1993}) through this transformation to a convergence theory on $\mathbb R$ (see Theorem~\ref{Theorem:Convergence}). The singular nature of the change of variables between these spaces is the main complication for the analysis. We make use of Sobolev spaces on the periodic interval $\mathbb T$ and $\mathbb R$.  We produce sufficient conditions for spectral convergence (faster than any polynomial) in various function spaces.  We do not address geometric rates of convergence, although that would be a natural extension of what we describe here.

The second goal of the current work is to use a similar approach for interpolatory projections.  The Dirichlet kernel can be used for interpolation on $\mathbb T$.  Furthermore, its $L^1$ norm produces the famous Lebesgue constant \cite{Galeev1983,Rivlin} and a bound on the trigonometric interpolation operator when acting on continuous, periodic functions.  Here we use estimates of the $L^p$ norm of the Dirichlet kernel and its first derivative to estimate the norm of the rational interpolation operator on $\mathbb R$.  We make use of results from \cite{Galeev1983} and modify them for our purposes in Appendix~\ref{Appendix:LpEstimates}.

The last goal of this paper is to use the estimates of Theorem~\ref{Theorem:Convergence} to present error bounds for formulas that appeared in \cite{Trogdon2013}.  In particular, we obtain Sobolev convergence to the boundary values of Cauchy operators and weighted $L^2$ convergence to the Fourier transform of smooth functions. The method for the Fourier transform is also shown to be asymptotic: the absolute error of the method decreases for increasing modulus of the wave number.

The paper is organized as follows.  We present background material on function spaces, convergence of periodic interpolants and the mechanics of interpolation on $\mathbb R$ in Sections~\ref{Section:T} and \ref{Section:R}. While these results are not new, we present them in order to keep the current work as self-contained and educational as possible.  In Section~\ref{Section:Convergence} we state and prove our main convergence theorem.  This is followed by Section~\ref{Section:Operator} which contains the main result concerning the norms of the interpolation operator.  Methods for both the oscillatory Cauchy integral (Section~\ref{Section:OscCauchy}) and Fourier transform (Section~\ref{Section:Fourier}) are then described with the relevant error bounds.  Appendices containing estimates on the Dirichlet kernel, sufficient conditions for convergence and numerical examples are included.  Throughout this manuscript we reserve the letter $C$ with and without subscripts for a generic constant that may vary from line to line.  Subscripts are used to denote any dependencies of the constant.

\begin{remark}
We concentrate on sufficient conditions in the current work.  We do not prove our estimates are sharp, just good enough to ensure convergence in many cases encountered in practice where one use these methods.
\end{remark}

\section{Function spaces and Interpolation on $\mathbb T$}\label{Section:T}

Some aspects of  periodic Sobolev spaces on the interval $\mathbb T$ are now reviewed.  We then review the convergence theory for trigonometric interpolants on the interval $\mathbb T$.  In this case, for the reader's benefit, we develop the theory from first principles.  For an interval $\mathbb I \subset \mathbb R$, define
\begin{align*}
L^p(\mathbb I) = \left\{ f : \mathbb I \goto \mathbb C, ~~ \text{measurable} ~:~ \int_{\mathbb I} |f(x)|^p dx < \infty\right\},
\end{align*}
with the norm
\begin{align*}
\|f\|_{L^p(\mathbb I)} = \left(\int_{\mathbb I} |f(x)|^p dx  \right)^{1/p}.
\end{align*}

\subsection{Periodic Sobolev spaces}

For $F \in L^2(\mathbb T)$, we use the notation
\begin{align}\label{FourierCoef}
\hat F_k = \frac{1}{2\pi} \int_0^{2\pi} e^{-ik\theta} F(\theta) \D\theta.
\end{align}
\begin{definition}
The periodic Sobolev space of order $s$ is defined by
\begin{align*}
H^s(\mathbb T) = \left\{ F \in L^2(\mathbb T) ~:~ \sum_{k=-\infty}^\infty (1+|k|)^{2s} |\hat F_k|^2 < \infty\right\},
\end{align*}
with the norm
\begin{align*}
\|F\|_{H^s(\mathbb T)} = \left( \sum_{k=-\infty}^\infty (1+|k|)^{2s} |\hat F_k|^2 \right)^{1/2}.
\end{align*}
\end{definition}

Note that $H^0(\mathbb T) = L^2(\mathbb T)$.  It  is well known that $H^s(\mathbb T)$ is a Hilbert space and it is a natural space (almost by definition) in which to study both Fourier series and trigonometric interpolants.  We use $\mathcal D$ to denote the (weak and strong) differentiation operator.  The usual prime notation is used for (strong) derivatives when convenient.

\begin{theorem}[\cite{atkinson}]\label{Theorem:Diff}
For $s \in \mathbb N$ and $F \in H^s(\mathbb T)$, $\mathcal D^jF(\theta)$ exists a.e. for $j = 1, \ldots, s$ and
\begin{align*}
\|\mathcal D^j F\|_{L^2(\mathbb T)} < \infty,\quad j = 1, \ldots, s.
\end{align*}
Furthermore,
\begin{align*}
|\|F\||_s = \left( \sum_{j=0}^s \|\mathcal D^j F\|^2_{L^2(\mathbb T)} \right)^{1/2},
\end{align*}
is a norm on $H^s(\mathbb T)$ and it is equivalent to the $H^s(\mathbb T)$ norm.
\end{theorem}

We also relate the $H^s(\mathbb T)$ space to spaces of differentiable functions.

\begin{definition}
The space of differentiable functions of order $r$ is defined by 
\begin{align*}
C^r_p(\mathbb T) = \left\{ F: \mathbb T \goto \mathbb C ~:~ F(\theta) = F( \theta + 2\pi), ~~ \mathcal D^r F(\theta) \text{ is continuous on } \mathbb T\right\},
\end{align*}
with norm
\begin{align*}
\|F\|_{C^r_p(\mathbb T)} = \sum_{j=0}^r \|\mathcal D^j F\|_u, ~~ \|F\|_u = \sup_{\theta \in \mathbb T} |F(\theta)|.
\end{align*}
\end{definition}

We require a common embedding result.
\begin{proposition}\label{Prop:PeriodicEmbed}
For $F \in H^s(\mathbb T)$ and $s > r + 1/2$
\begin{align*}
\|F\|_{C_p^r(\mathbb T)} \leq C_{s,r} \|F\|_{H^s(\mathbb T)}.
\end{align*}
\begin{proof}
Since $r \geq 0$, we have that $ s > 1/2$.  Therefore
\begin{align*}
|F(\theta)| \leq \sum_k |e^{ik\theta} \hat F_k| \leq \left(\sum_k (1+|k|)^{-2s} \right)^{1/2} \left(\sum_k |\hat F_k|^2(1+|k|)^{2s}\right)^{1/2},
\end{align*}
by the Cauchy-Schwarz inequality.  This shows two things.  First, because $s > 1/2$, $F(\theta)$ is the uniform limit of continuous functions and is therefore continuous.  Second, $\|F\|_u \leq C_s \|F\|_{H^s(\mathbb T)}$.  More generally, for $j \leq r$
\begin{align*}
|\mathcal D^j F(\theta)| &\leq \sum_k |(ik)^j e^{ik\theta} \hat F_k| \leq 
\sum_k (1+|k|)^j |\hat F_k|\\
&\leq \left(\sum_k (1+|k|)^{-2(s-j)} \right)^{1/2} \left(\sum_k |\hat F_k|^2(1+|k|)^{2s}\right)^{1/2}.
\end{align*}
The first sum in the last expression converges because $s - r > 1/2$ ($j \leq r$). Taking a supremum we find
\begin{align*}
\|\mathcal D^j F\|_u \leq  C_{s,j} \|F\|_{H^s(\mathbb T)}. 
\end{align*}
This proves the result.
\end{proof}
\end{proposition} 

\begin{remark}
We have ignored some technicalities in the proof of the previous result.  Because an $H^s(\mathbb T)$ function may take arbitrary values on a Lebesgue set of measure zero, we need to make sense of $\|F\|_u$.  It suffices to take the representation of $F \in H^s(\mathbb T)$, $s > 1/2$ defined by its Fourier series because we are guaranteed that this is a continuous function.
\end{remark}

\subsection{Trigonometric interpolation}\label{Section:Interp}

We now discuss the construction of the series representations for trigonometric interpolants of a continuous function $F$. For $n \in \mathbb N$, define $\theta_j = 2 \pi j / n$ for $j = 0, \ldots, n-1$ and two positive integers
\begin{align*}
n_+ = \lfloor n/2 \rfloor, \quad n_- = \lfloor (n-1)/2 \rfloor.
\end{align*}
Note that $n_+ + n_- + 1 = n$ regardless of whether $n$ is even or odd.

\begin{definition}
The discrete Fourier transform of order $n$ of a continuous function $F$ is the mapping
\begin{align*}
\mathcal F_n F &= ( \tilde F_0, \tilde F_1, \ldots, \tilde F_n)^T,\\
\tilde F_k &= \frac{1}{n} \sum_{j=0}^{n-1} e^{-ik\theta_j} F(\theta_j).
\end{align*}
\end{definition}

This is nothing more than the trapezoidal rule applied to \eqref{FourierCoef}.  Note that the dependence on $n$ is implicit in the $\tilde F_k$ notation.  These coefficients, unlike true Fourier coefficients, depend on $n$.  This formula produces the coefficients for the interpolant.

\begin{proposition}\label{Prop:Interpolation}
The function
\begin{align*}
\mathcal I_n F(\theta) = \sum_{k=-n_-}^{n_+} e^{ik\theta} \tilde F_k,
\end{align*}
satisfies $\mathcal I_n F(\theta_j) = F(\theta_j)$.
\begin{proof}
By direct computation
\begin{align*}
\mathcal I_n F(\theta_j ) = \sum_{k=-n_-}^{n_+} e^{ik\theta_j} \left( \frac{1}{n} \sum_{\ell=0}^{n-1} e^{-ik\theta_\ell} F(\theta_\ell)  \right) = \frac{1}{n} \sum_{\ell=0}^{n-1} F(\theta_\ell) \sum_{k=-n_-}^{n_+} e^{ik(\theta_j-\theta_\ell)}.
\end{align*}
Therefore we must compute
\begin{align*}
\sum_{k=-n_-}^{n_+} e^{ik(\theta_j-\theta_\ell)} = e^{- ikn_-(\theta_j-\theta_\ell)} \sum_{k=0}^{n-1} e^{ik(\theta_j-\theta_\ell)},
\end{align*}
using $n_++n_- = n-1$.  Then it is clear that
\begin{align*}
\sum_{k=0}^{n-1} e^{ik(\theta_j-\theta_\ell)} = n \delta_{j\ell},
\end{align*}
where $\delta_{j\ell}$ is the usual Kronecker delta.  This is easily seen using the formula for the partial sum of a geometric series.  This proves the result.
\end{proof}
\end{proposition}

From these results it is straightforward to construct the matrix (and its inverse!) that maps the vector $(F(\theta_1),\ldots,F(\theta_n))^T$ to $( \tilde F_0, \ldots, \tilde F_n)^T$ but we will not pursue this further other than to say that this mapping is efficiently computed with the FFT.  We have used suggestive notation above: $\mathcal I_n$ is used to denote the projection operator that maps $F \in C_p^0(\mathbb T)$ to its (unique) interpolant $\mathcal I_n F$.

Reviewing the proof of Proposition~\ref{Prop:Interpolation} we see that there is an alternate expression for $\mathcal I_n F(\theta)$:
\begin{align}
\mathcal I_nF(\theta) &= \sum_{\ell= 0}^{n-1} F(\theta_\ell) \frac{D_n(\theta-\theta_\ell)}{n},\label{OpForm}\\
D_n(\theta) &= \sum_{k=-n_-}^{n_+} e^{ik\theta}.\label{Dirichlet}
\end{align}
The function $D_n$ is referred to in the literature as the Dirichlet kernel and it plays a central role in Section~\ref{Section:Operator}.  We now present the result of Kress and Sloan \cite{Kress1993} for the convergence of trigonometric interpolants.

\begin{theorem}[\cite{Kress1993}]\label{Theorem:DFT-converge}
For $F \in H^s(\mathbb T)$, $s > 1/2$ and $0 \leq t \leq s$,
\begin{align*}
\|\mathcal I_n F-F\|_{H^t(\mathbb T)} \leq C_{t,s} n^{t-s} \|F\|_{H^s(\mathbb T)}.
\end{align*}
\begin{proof}
We begin with the expression
\begin{align*}
F(\theta) = \sum_k e^{ik\theta} \hat F_k.
\end{align*}
Let $e_k(\theta) = e^{ik\theta}$ and consider the interpolation of these exponentials. For $-n_- \leq k \leq n_+$ it is clear that $\mathcal I_n e_k = e_k$. For other values of $k$, write $k = j + mn$ for $m \in \mathbb N \setminus\{0\}$ and $-n_- \leq j \leq n_+$.  Then
\begin{align*}
e_{j+mn}(\theta_\ell) = \exp(2\pi i/n \ell(j+mn)) = \exp (2\pi i\ell m + 2 \pi i j\ell/n) = e_j(\theta_\ell).
\end{align*}
This is the usual aliasing relation.  From this we conclude $\mathcal I_n e_{j+mn} = e_j$. We now consider applying the interpolation operator to the Fourier series expression for $F$:
\begin{align*}
\mathcal I_nF(\theta) = \sum_{k=-n_-}^{n_+}  e^{ik\theta} \hat F_k + {\sum_{m=-\infty}^\infty}' \left( \sum_{k=-n_-}^{n_+}  e^{ik\theta} \hat F_{k+mn} \right),
\end{align*}
were the $'$ indicates that the $m = 0$ term is omitted in the sum.  We find
\begin{align*}
\mathcal I_nF(\theta)-F(\theta) = {\sum_{m=-\infty}^\infty}' \left( \sum_{k=-n_-}^{n_+}  e^{ik\theta} \hat F_{k+mn} \right) - \left( \sum_{k=-\infty}^{-n_--1} + \sum_{k=n_++1}^{\infty} \right)e^{ik\theta} \hat F_k.
\end{align*}
We estimate the $H^t(\mathbb R)$ norms of the two terms individually.  First,
\begin{align*}
S_+(\theta) &= \sum_{k=n_++1}^{\infty} e^{ik\theta} \hat F_k,\\
\|S_+\|_{H^t(\mathbb R)}^2 &= \sum_{k=n_++1}^{\infty} (1+|k|)^{2t} |\hat F_k|^2 = \sum_{k=n_++1}^{\infty} (1+|k|)^{2(t-s)} (1+|k|)^{2s}|\hat F_k|^2\\
&\leq (2 + n_+)^{2(t-s)} \|F\|^2_{H^s(\mathbb T)}.
\end{align*}
The same estimate holds for $S_-(\theta) = \sum_{k=-\infty}^{-n_--1} e^{ik\theta} \hat F_k$ with $n_+$ replaced with $n_-$.  It remains to estimate (after switching summations)
\begin{align*}
S_0(\theta) = \sum_{k=-n_-}^{n_+}e^{ik\theta} {\sum_{m=-\infty}^\infty}'    \hat F_{k+mn}
\end{align*}
so that
\begin{align*}
\|S_0\|_{H^t(\mathbb T)}^2 = \sum_{k=-n_-}^{n_+} (1+|k|)^{2t}\left|{\sum_{m=-\infty}^\infty}'    \hat F_{k+mn}\right|^2.
\end{align*}
It follows that
\begin{align*}
\left|{\sum_{m=-\infty}^\infty}'    \hat F_{k+mn}\right| &\leq {\sum_{m=-\infty}^\infty}'    (1+|k+mn|)^{-s} (1+|k+mn|)^s |\hat F_{k+mn}| \\
& \leq \left({\sum_{m=-\infty}^\infty}' (1+|k+mn|)^{-2s}\right)^{1/2} \left({\sum_{m=-\infty}^\infty}' (1+|k+mn|)^{2s}|\hat F_{k+mn}|^2\right)^{1/2}\\
&\leq n^{-s}\left({\sum_{m=-\infty}^\infty}' (|k/n+m|)^{-2s}\right)^{1/2} \left({\sum_{m=-\infty}^\infty}' (1+|k+mn|)^{2s}|\hat F_{k+mn}|^2\right)^{1/2}.
\end{align*}
Observe that for $t \in [-1/2,1/2]$ function
\begin{align*}
\left({\sum_{m=-\infty}^\infty}' (|t+m|)^{-2s}\right)^{1/2}, ~~ s > 1/2,
\end{align*}
is continuous and is therefore bounded uniformly by a constant $c_s$. Then
\begin{align*}
\|S_0\|_{H^t(\mathbb T)}^2 \leq c_s^2n^{-2s} \sum_{k=-n_-}^{n_+} (1+|k|)^{2t} {\sum_{m=-\infty}^\infty}' (1+|k+mn|)^{2s}|\hat F_{k+mn}|^2 \leq  c_s^2n^{-2s}(1+n_+)^{2t} \|F\|_{H^s(\mathbb T)}.
\end{align*}
Combining the estimates for $S_\pm$ and $S_0$ proves the result.
\end{proof}
\end{theorem}

\section{Function Spaces and Interpolation on $\mathbb R$}\label{Section:R}

We mirror the previous section and present results for spaces of functions defined on $\mathbb R$.  We use a M\"obius transformation to construct a rational interpolant of a continuous function on $\mathbb R$.

\subsection{Function spaces on the line}

Our major goal is the rational approximation of functions defined on $\mathbb R$ and we introduce the relevant function spaces.  The Fourier transform (for $f \in L^2(\mathbb R)$) is defined by
\begin{align*}
\mathcal F f(k) &= \int_{\mathbb R} e^{-ikx} f(x) dx,\\
f(x) &= \frac{1}{2 \pi} \int_{\mathbb R} e^{ikx} \mathcal F f(k) dk.
\end{align*}
We present a series of results about these Sobolev spaces.  For the sake of brevity, we do not prove these here. A general reference is \cite{Folland}.

\begin{definition}\label{Def:SobolevLine}
The Sobolev space on the line of order $s$ is defined by
\begin{align*}
H^s(\mathbb R) = \left\{ f \in L^2(\mathbb R) ~:~ \int_{\mathbb R} (1+|k|)^{2s} |\mathcal F f(k)|^2 dk < \infty \right\},
\end{align*}
with the norm
\begin{align*}
\|f\|_{H^s(\mathbb R)} = \left(\int_{\mathbb R} (1+|k|)^{2s} |\mathcal F f(k)|^2 dk\right)^{1/2}.
\end{align*}
\end{definition}

\begin{theorem}[\cite{Folland}]
Theorem \ref{Theorem:Diff} holds with $\mathbb T$ replaced with $\mathbb R$.
\end{theorem}

\begin{definition}
The space of differentiable functions of order $r$ that decay at infinity is defined by
\begin{align*}
C_0^r(\mathbb R) = \left\{ f : \mathbb R \goto \mathbb C ~:~ \mathcal D^r f(x) \text{ is continuous on } \mathbb R, \lim_{|x| \goto \infty} \mathcal D^jf(x) = 0, ~~ j=0,1,\ldots,r\right\}.
\end{align*}
with norm
\begin{align*}
\|f\|_{C^r_0(\mathbb R)} = \sum_{j=0}^r \|\mathcal D^j f\|_u, ~~ \|f\|_u = \sup_{x \in \mathbb R} |f(x)|.
\end{align*}
\end{definition}

\begin{theorem}[Sobolev Embedding,\cite{Folland}]
For $f \in H^s(\mathbb R)$ and $ s > r + 1/2$ then $f \in C_0^r(\mathbb R)$ and
\begin{align*}
\|f\|_{C_0^r(\mathbb R)} \leq C_{r,s} \|f\|_{H^s(\mathbb R)}.
\end{align*}
\end{theorem}

\subsection{Practical Rational Approximation}

In this section, a method for the rational approximation of functions on $\mathbb R$ is discussed. The fundamental tool is the FFT that was discussed in the previous section.  Two references for this method are \cite{SOHilbertTransform,Trogdon2013} although we follow \cite{Trogdon2013} closely.  Define a one-parameter family of M\"obius transformations
\begin{align*}
M_\beta(z) = \frac{z - i \beta}{z+i\beta}, ~~ M^{-1}_\beta(z) = \frac{\beta}{i} \frac{z+1}{z-1},~~ \beta >0.
\end{align*}
Each of these transformations maps the real axis to the unit circle. Assume $f$ is a smooth and rapidly decaying function on $\mathbb R$.  Then $f$ is mapped to a smooth function on $[0,2\pi]$ by $F(\theta) = f(M_\beta^{-1}(e^{i\theta}))$ (see Proposition~\ref{Prop:Decay}).  Thus, the FFT may be applied to $F(\theta)$ to obtain a sequence $\mathcal I_nF(\theta)$ of rapidly converging interpolants. The transformation $x = T(\theta) = M_\beta^{-1}(e^{i\theta})$ is inverted:
\begin{align*}
\mathcal R_nf(x) = \mathcal I_n F(T^{-1}(x)),
\end{align*}
is a rational approximation of $f$.  We examine this expansion more closely.

From Section~\ref{Section:Interp} we have
\begin{align*}
\mathcal I_n F(\theta) &= \sum_{k=-n_-}^{n_+} e^{ik\theta} \tilde F_k,\\
\tilde F_k & = \frac{1}{n} \sum_{\ell = 0}^{n-1} e^{ik\theta_\ell} f(M_\beta^{-1}(e^{i\theta_\ell})).
\end{align*}
Then, in mapping to the real axis we find
\begin{align*}
\mathcal R_nf(x) = \sum_{k=-n_-}^{n_+}\tilde F_k M_\beta^k(x).
\end{align*}

\begin{remark}
Even though we express $\tilde F_k$ as a sum, note that in practice it should be computed with the FFT for efficiency.
\end{remark}

The behavior of $\mathcal R_nf$ at $\infty$ is important.  It is clear that $\lim_{\theta \goto 0^+} M_\beta^{-1}(e^{i\theta}) = +\infty$ so that $\mathcal I_nF(0) = 0$.  This implies that $\sum_{k=-n_-}^{n_+} \tilde F_k = 0$ and
\begin{align*}
\mathcal R_nf(x) = \sum_{k=-n_-}^{n_+}\tilde F_k (M_\beta^k(x)-1).
\end{align*}
In following with \cite{Trogdon2013}, we drop $\beta$ dependence and define $R_k(x) = M_\beta^k(x)-1$.  In summary, we have designed a method for the interpolation  of a function in the basis $\{R_k\}_{k=-\infty}^\infty$.  Indeed, this is a basis of $L^2(\mathbb R)$ once $R_0(x) = 0$ is removed \cite{Trogdon2013}.  We devote the entire next section to the study of convergence.

\section{Convergence}\label{Section:Convergence}

We discuss various convergence properties of the sequence $\{\mathcal R_nf\}_{n >1}$ depending on the regularity of $f$.  As is natural, all properties are derived from the convergence of the discrete Fourier transform.  Throughout this section, and the remainder of the manuscript, we associate $f$ and $F$ by the change of variables $F(\theta) = f(T(\theta))$, $T(\theta) = M_\beta^{-1}(e^{i\theta})$.  We summarize the results of this section in the following theorem.

\begin{theorem}\label{Theorem:Convergence}
Assume $F \in H^s(\mathbb R)$, $s > 1/2$ and $f$ is in appropriate function spaces to make the following norms finite.  Then:
\begin{itemize}
\item $\|\mathcal R_n f- f\|_{C_0^r(\mathbb R)} \leq C_{r,s} n^{1/2+r-s} \|F\|_{H^s(\mathbb T)}$, ~~ $r < s + 1/2$, and
\item $\|\mathcal R_n f- f\|_{H^t(\mathbb R)} \leq [C_{\epsilon,s} n^{1+\epsilon-s} + C_{s,t} n^{t-s} ]\|F\|_{H^s(\mathbb T)}, ~~ \epsilon > 0, ~~t < s$, ~$s > 1 + \epsilon$,
\item $\displaystyle \left| \dashint_{\mathbb R} (\mathcal R_nf(x) - f(x)) dx \right| \leq C_{\epsilon,s} n^{3/2+\epsilon-s}, ~~ \epsilon > 0$,~ $s > 3/2 + \epsilon$, and 
\item $\|\mathcal D^j [\mathcal R_n f- f]\|_{L^1(\mathbb R)} \leq C_{j,s} n^{j-s} \|F\|_{H^s(\mathbb T)}$, $j > 0$.
\end{itemize}
Here constants $C_{\epsilon,s}$ (which may differ in each line) are unbounded as $\epsilon \goto 0^+$.
\end{theorem}

Before we prove each of the estimates we need a lemma.
\begin{lemma}\label{Lemma:Vanish}
For $f \in C_0^r(\mathbb T)$, $f^{(j)}(T(\theta)) = \sum_{\ell =1}^j F^{(\ell)}(\theta) p_\ell(\theta)$ where $p_\ell(\theta)$ is bounded and vanishes to at least second order at $\theta = 0,2\pi$.
\begin{proof}
First, observe that
\begin{align*}
f'(x) &= F'(T^{-1}(x)) \frac{d}{dx} T^{-1}(x),\\
f'(T(\theta)) &= F'(\theta) [T'(\theta)]^{-1}.
\end{align*}
The general case is seen by showing that $ \frac{d}{dx} T^{-1}(x)$ and all its derivatives decay are $\bigo(x^{-2})$ as $|x| \goto \infty$.
\end{proof}
\end{lemma}

We prove each piece of the theorem in a subsection.

\subsection{Uniform convergence}

From Theorem~\ref{Theorem:DFT-converge} we have $\|\mathcal I_n F-F\|_{H^t(\mathbb T)} \leq C_{t,s} n^{t-s} \|F\|_{H^s(\mathbb T)}$ and combining this with Proposition~\ref{Prop:PeriodicEmbed} we find
\begin{align*}
\|\mathcal I_n F-F\|_{C_p^0(\mathbb T)} = \|\mathcal R_n f-f\|_{C_0^0(\mathbb R)} \leq C_{s} n^{1/2-s} \|F\|_{H^s(\mathbb T)}.
\end{align*}
Therefore, we easily obtain uniform convergence provided the mapped function $F$ is smooth. Furthermore, using Lemma~\ref{Lemma:Vanish}, Proposition~\ref{Prop:PeriodicEmbed} and Theorem~\ref{Theorem:DFT-converge} we find
\begin{align*}
\|\mathcal R_n f-f\|_{C_r^0(\mathbb R)} \leq C_r \|\mathcal I_nF-F\|_{C_p^r(\mathbb T)} \leq C_{r,s} n^{1/2+r-s} \|F\|_{H^s(\mathbb T)}.
\end{align*}

\subsection{$L^2(\mathbb R)$ convergence}

Demonstrating the convergence of the approximation in $L^2(\mathbb R)$ is a more delicate procedure.  We directly consider
\begin{align*}
\|\mathcal R_nf-f\|^2_{L^2(\mathbb R)} =\int_{\mathbb R} |\mathcal R_nf(x) - f(x)|^2 dx = \int_0^{2\pi}  |\mathcal R_nf(T(\theta)) - f(T(\theta))|^2 |dT(\theta)|.
\end{align*}
We find $T'(\theta) = \frac{2}{\beta} \frac{e^{i\theta}}{(e^{i\theta}-1)^2}$ so that
\begin{align*}
\|\mathcal R_nf-f\|^2_{L^2(\mathbb R)}= \frac{2}{\beta}\int_0^{2\pi}  |\mathcal I_nF(\theta) - F(\theta)|^2 \frac{d\theta}{|e^{i\theta}-1|^2}.
\end{align*}
The unbounded nature of the change of variables makes it clear that we must require the convergence of a derivative.  We break the integral up into two pieces.  For the first piece we write $H_n(\theta) = \mathcal I_nF(\theta) - F(\theta)$ while noting that $H_n(0) = 0$
\begin{align*}
\frac{2}{\beta}\int_0^{\pi}  \left|\int_0^\theta H_n'(\theta') d\theta' \right|^2 \frac{d\theta}{|e^{i\theta}-1|^2} \leq \left( \int_{0}^{2\pi} |H_n'(\theta')|^{p} d\theta'\right)^{2/p} \frac{2}{\beta}\int_0^{\pi} \frac{|\theta|^{2/q}}{|e^{i\theta}-1|^2} d\theta,
\end{align*}
for $1/p + 1/q = 1$. A similar estimate holds for the integral from $\pi$ to $2\pi$.  It is clear that $q< 2$ is required for the integral to converge.  It remains to express the integral involving $H_n'$ in terms of something known.  A well-known fact is that if $\|G\|_{L^2(\mathbb T)} \leq c_2$ and $\|G\|_{C_p^0(\mathbb T)} \leq c_u$ then for $2 \leq p < \infty$
\begin{align*}
\|G\|_{L^p(\mathbb T)} \leq c_u^{1-2/p}c_2^{2/p}.
\end{align*}
We find
\begin{align*}
\|H_n'\|_{C_p^0(\mathbb T)} \leq C_s n^{3/2-s} \|F\|_{H^s(\mathbb T)}
\end{align*}
from Theorem~\ref{Theorem:DFT-converge} and Proposition~\ref{Prop:PeriodicEmbed}.  Also,
\begin{align*}
\|H_n'\|_{L^2(\mathbb T)} \leq C_s n^{1-s} \|F\|_{H^s(\mathbb T)},
\end{align*}
from Theorem~\ref{Theorem:DFT-converge}.  Then
\begin{align*}
\|H_n'\|_{L^p(\mathbb T)} \leq C_{p,s} n^{3/2-1/p-s} \|F\|_{H^s(\mathbb T)},~~ p > 2,
\end{align*}
which results in
\begin{align*}
\|\mathcal R_nf - f\|_{L^2(\mathbb R)} \leq C_{\epsilon,s} n^{1+\epsilon-s} \|F\|_{H^s(\mathbb T)}, ~~ \epsilon > 0.
\end{align*}

\subsection{$H^t(\mathbb R)$ convergence}

We use Lemma~\ref{Lemma:Vanish} and directly compute,
\begin{align}\label{derivative-series}
\|f^{(j)}\|_{L^2(\mathbb R)} \leq \sum_{\ell =1}^j \left( \int_{0}^{2\pi} |F^{(\ell)}(\theta)|^2 |p_\ell(\theta)|^2 |T'(\theta)|\D\theta \right)^{1/2} \leq C_j \|F\|_{H^j(\mathbb T)}.
\end{align}
Therefore
\begin{align*}
\|f\|_{H^t(\mathbb R)} \leq \|f\|_{L^2(\mathbb R)} + C_t \|F\|_{H^t(\mathbb T)}.
\end{align*}
Replacing $f$ with $\mathcal R_nf-f$ we have
\begin{align*}
\|\mathcal R_n f-f\|_{H^t(\mathbb R)} \leq [C_{\epsilon,s} n^{1+\epsilon-s} + C_{s,t} n^{t-s} ]\|F\|_{H^s(\mathbb T)}, ~~ \epsilon > 0.
\end{align*}

\begin{remark}
This bound seems to indicate that convergence of the function in $L^2(\mathbb R)$ requires more smoothness than convergence of the first derivative in $L^2(\mathbb R)$. This is an artifact of using the smoothness of the mapped function $F$ to measure the convergence rate.
\end{remark}

\subsection{Convergence of the integral}\label{Section:ConvIntegral}

Showing convergence of the integral of $\mathcal R_nf$ to that of $f$ is even more delicate than $L^2(\mathbb R)$ convergence.  The main reason for this is that generically $\mathcal R_nf \not\in L^1(\mathbb R)$ although it has a convergent principal-value integral.   Therefore the quantity of study is
\begin{align*}
S_n(f) = \left| \dashint_{\mathbb R} (\mathcal R_nf(x) - f(x)) dx \right|,
\end{align*}
where
\begin{align*}
\dashint f(x) dx = \lim_{R \goto \infty} \int_{-R}^R f(x) dx,
\end{align*}
if the limit exists.  Let $H_n$ be as in the previous section and let $s_n = H_n'(0)$, assuming that $F$ has at least one continuous derivative.  Define $\tilde H_n(\theta) = H_n(\theta) + is_n(e^{i\theta}-1)$ and note that  $s_n$ is chosen so that $\tilde H_n(0) = \tilde H_n'(0) = 0$.  The derivative also vanishes for $\theta =2\pi$. Turning back to $S_n(f)$ we have
\begin{align}\label{fixed-pv} 
S_n(f) &= \left| \int_{\mathbb R} (\mathcal R_nf(x) + i s_n R_1(x) - f(x)) dx  - i s_n \dashint R_1(x) dx \right|\\
&\leq \left|\int_{\mathbb R} (\mathcal R_nf(x) + i s_n R_1(x) - f(x)) dx  \right| + 2\pi \beta |s_n|.
\end{align}
Changing variables on the integral, we have
\begin{align*}
\left|\int_{\mathbb R} (\mathcal R_nf(x) + i s_n R_1(x) - f(x)) dx  \right| \leq \int_0^{2\pi} \frac{|\tilde H_n(\theta)|}{|e^{i\theta}-1|^2} d\theta.
\end{align*}
A straightforward estimate produces
\begin{align*}
|\tilde H_n(\theta)| \leq  \int_{0}^{\theta}\int_{0}^{\theta'} |\tilde H_n''(\theta'')| \D\theta''\D\theta' \leq \|\tilde H_n''\|_{L^p(\mathbb T)} \frac{\theta^{1+1/q}}{1 + 1/q},
\end{align*}
for $1/p+1/q = 1$. If we require that $p > 1$ ($q < \infty$) then 
\begin{align*}
\int_0^{2\pi} \frac{|\tilde H_n(\theta)|}{|e^{i\theta}-1|^2} d\theta \leq C_p \|\tilde H_n''\|_{L^p(\mathbb T)}.
\end{align*}
Therefore
\begin{align*}
\|\tilde H_n''\|_{L^p(\mathbb T)} \leq \|H_n''\|_{L^p(\mathbb T)} + C |s_n|.
\end{align*}
Similar estimates to those in the previous section produce for $p > 1$
\begin{align*}
\|H_n''\|_{L^p(\mathbb T)} &\leq C_{p,s} n^{5/2-1/p-s} \|F\|_{H^s(\mathbb T)},\\
|s_n| &\leq C_s n^{3/2-s} \|F\|_{H^s(\mathbb T)}.
\end{align*}
Therefore
\begin{align*}
S_n(f) \leq C_{\epsilon,s} n^{3/2+\epsilon-s}, ~~ \epsilon > 0.
\end{align*}

To consider derivatives, we note that the principal-value integral is no longer needed.  For $f \in H^j(\mathbb R)$, $\mathcal D^jf \in L^1(\mathbb R)$, $j \geq 1$ (see \eqref{derivative-series})
\begin{align*}
\|\mathcal D^j f\|_{L^1(\mathbb R)} \leq  C_j \sum_{\ell=0}^j \|\mathcal D^\ell F\|_{L^1(\mathbb T)}.  
\end{align*}
Replacing $f$ with $\mathcal R_n f -f$  and noting that the $L^2(\mathbb T)$ norm dominates the $L^1(\mathbb T)$ norm we find
\begin{align*}
\|D^j[f-\mathcal R_nf]\|_{L^1(\mathbb R)} \leq  C_j\|F-\mathcal I_n F\|_{H^j(\mathbb T)} \leq C_{j,s} n^{j-s} \|F\|_{H^s(\mathbb T)}.
\end{align*}

\section{An Interpolation Operator}\label{Section:Operator}

In this section we prove the following result concerning the norm of $\mathcal R_n$.

\begin{theorem}\label{Theorem:Operator}
The interpolation projection operator $\mathcal R_n$ satisfies the following estimates for $n > 1$
\begin{itemize}
\item $\|\mathcal R_n\|_{C_0^0(\mathbb R) \goto C_0^0(\mathbb R)} \leq C \log n $,
\item $\|\mathcal R_n\|_{C_0^0(\mathbb R) \goto L^2(\mathbb R)} \leq C n\log n$,
\item $\|\mathcal R_n\|_{C_0^0(\mathbb R) \goto H^1(\mathbb R)} \leq C n^{3/2}$, and hence
\item $\|\mathcal R_n\|_{H^1(\mathbb R) \goto H^1(\mathbb R)} \leq C n^{3/2}$.
\end{itemize}
\end{theorem}

The main result we need here to prove this theorem is the following estimates of the Dirichlet kernel $D_n$ (see \eqref{Dirichlet}).  This theorem, as stated, is a special case of the general results of \cite{Galeev1983}.
\begin{theorem}[\cite{Galeev1983}]\label{Theorem:Galeev}
Let $\alpha \in \mathbb N$.  Then for $n >1$
\begin{align*}
\|D_n^{(\alpha)}\|_{L^p(\mathbb T)} \leq \begin{choices} C_{\alpha,p} n^{\alpha+1-1/p}, \when \alpha + 1 - 1/p > 0,\\
\\
C \log n, \when \alpha = 0, ~p=1,\\
\end{choices}
\end{align*}
\end{theorem}
\noindent This theorem is proved for $1< p< \infty$ and $\alpha = 0,1$ in Appendix~\ref{Appendix:LpEstimates}. As in the previous section we prove each piece of the theorem in a subsection.

\subsection{Uniform operator norm}
This estimate derived here is just the usual Lebesgue constant for interpolation \cite{atkinson}.  As before, uniform bounds on functions translate directly between the real axis and $\mathbb T$.  It is clear that when considering \eqref{OpForm}
\begin{align*}
\|\mathcal R_nf\|_{C_0^0(\mathbb R)} = \|\mathcal I_nF\|_{C_p^0(\mathbb T)} =  \frac{1}{n} \|F\|_{C_p^0(\mathbb T)} \|D_n\|_{L^1(\mathbb T)} \leq C \log n\|F\|_{C_p^0(\mathbb T)}.
\end{align*}
Therefore 
\begin{align*}
\|\mathcal R_n\|_{C_0^0(\mathbb R) \goto C_0^0(\mathbb R)} \leq C \log n.
\end{align*}

\subsection{$L^2(\mathbb R)$ and $H^1(\mathbb R)$ operator norms}

Estimates on the $L^2(\mathbb R)$ and $H^1(\mathbb R)$ operator norms require more care.  Because $\mathcal I_nF(\theta) = 0$ in the case that $f \in H^1(\mathbb R)$ or $f \in C_0^0(\mathbb R)$ we write
\begin{align*}
\mathcal I_nF(\theta) &= \sum_{\ell=0}^{n-1} F(\theta_\ell) \frac{1}{n} \left( D_n(\theta-\theta_\ell) - D_n(-\theta_\ell) \right) = \sum_{\ell=0}^{n-1} F(\theta_\ell) \frac{1}{n} \int_{0}^\theta D_n'(\theta'-\theta_\ell) d \theta'\\
&= \sum_{\ell=0}^{n-1} F(\theta_\ell) \frac{1}{n} \int_{2\pi}^\theta D_n'(\theta'-\theta_\ell) d \theta'.
\end{align*}

Using the change of variables $T$, we find
\begin{align}\label{l2norm}
\|\mathcal R_nf\|_{L^2(\mathbb R)} \leq \frac{2}{\beta} \|F\|_{C_p^0(\mathbb T)} \frac{1}{n} \sum_{\ell=0}^{n-1} \left[ \left( \int_{0}^{\pi} \frac{ \left|\int_{0}^{\theta} D_n'(\theta'-\theta_\ell) d\theta'\right|^2}{|e^{i\theta}-1|^2} d\theta \right)^{1/2} + \left( \int_{\pi}^{2\pi} \frac{ \left|\int_{2 \pi}^{\theta} D_n'(\theta'-\theta_\ell) d\theta'\right|^2}{|e^{i\theta}-1|^2} d\theta \right)^{1/2}\right].
\end{align}

The next step is the estimation of the integral
\begin{align*}
\int_{0}^{\pi} \frac{ \left|\int_{0}^{\theta} D_n'(\theta'-\theta_\ell) d\theta'\right|^2}{|e^{i\theta}-1|^2} d\theta.
\end{align*}
  We further break up this integral and consider
\begin{align*}
I_\ell = \int_{0}^{\theta_\ell/2} \frac{ \left|\int_{0}^{\theta} D_n'(\theta'-\theta_\ell) d\theta'\right|^2}{|e^{i\theta}-1|^2} d\theta\leq \int_0^{\theta_\ell/2} \frac{\theta^{2/q}}{|e^{i\theta}-1|^2} \left| \int_0^\theta |D_n'(\theta'-\theta_\ell)|^p d\theta' \right|^{2/p}  d \theta.
\end{align*}
It is clear that $q < 2$ is required for the integrability of the first factor, therefore $p  > 2$.  Also, $1/p+1/q = 1$.  The integration bounds are sufficient to ensure that the argument of $D_n'$ is bounded away from zero.  Furthermore, under these constraints, using that $D_n'$ is odd
\begin{align*}
\left|\int_0^\theta |D_n'(\theta'-\theta_\ell)|^p d\theta' \right|^{1/p} \leq \|D_n'\|_{L^p(\theta_\ell/2, \theta_\ell)} \leq \|D_n'\|_{L^p(\theta_\ell/2, 2\pi-\theta_\ell/2)},\\
\leq C_p(1 + \theta_\ell n_-/2) \left( \frac{n_-}{ 1 + \theta_\ell n_-/2} \right)^{2-1/p}.
\end{align*}
This estimate follows from \eqref{derivative-estimates} with $\epsilon/n_- = \theta_\ell/2 \leq \pi$.  Next, note that for $n \geq 1$, $\pi \ell/2 \leq \theta_\ell n_- = 2 \pi \ell \lfloor(n-1)/2\rfloor/n \leq \pi \ell$.  Then
\begin{align*}
\|D_n'\|_{L^p(\theta_\ell/2, 2\pi)} \leq C_p n_-^{2-1/p} \ell^{1/p-1}
\end{align*}
for a new constant $C_p$.  Next, it is clear that
\begin{align*}
\int_0^\theta \frac{\theta'^{2/q}}{|e^{i\theta'}-1|^2} d \theta' \leq C_q \theta^{2/q-1},
\end{align*}
so that
\begin{align}\label{inside-estimate}
I_\ell^{1/2} \leq C_{p,q}\theta_\ell^{1/q-1/2 } n_-^{2-1/p} \ell^{1/p-1} \leq C_{p,q} \pi^{1/q-1/2} n_-^{3/2} \ell^{-1/2}.
\end{align}
Here we used $1/p+1/q = 1$.

Next, we consider
\begin{align}
L_\ell &= \int_{\theta_\ell/2}^{\pi} \frac{ \left|\int_{0}^{\theta} D_n'(\theta'-\theta_\ell) d\theta'\right|^2}{|e^{i\theta}-1|^2} d\theta \leq \frac{4\pi }{\theta_\ell^2} \|D'_n\|_{L^1(\mathbb T)}^2 \leq C \frac{4\pi n^2 }{\theta_\ell^2} \leq C^2 \frac{n^4}{\ell^2},\notag\\
L_\ell^{1/2} &\leq C \frac{n^2}{\ell}. \label{outside-estimate}
\end{align}
by appealing to Theorem~\ref{Theorem:Galeev}.

Assembling \eqref{outside-estimate} and \eqref{inside-estimate} we find
\begin{align*}
\frac{1}{n} \sum_{\ell=0}^{n-1} \left( \int_{0}^{\pi} \frac{ \left|\int_{0}^{\theta} D_n'(\theta'-\theta_\ell) d\theta'\right|^2}{|e^{i\theta}-1|^2} d\theta \right) \leq \frac{C}{n} \sum_{\ell=0}^{n-1} [n^{3/2} \ell^{-1/2} + n^2 \ell^{-1}] \leq C'n (1 + \log n).
\end{align*}
It is clear the remaining integrals from $\pi$ to $2\pi$ in \eqref{l2norm} have a similar bound.  Therefore, we conclude that for $n \geq 1$
\begin{align*}
\|\mathcal R_n f\|_{L^2(\mathbb R)} \leq C n \log n \|f\|_{C_0^0(\mathbb R)}.
\end{align*}

\begin{remark}
The heuristic reason for the previous estimates is as follows.  We expect $D_n'(\theta-\theta_\ell)$ to be largest near $\theta_\ell$.  The integrand will be further amplified when $\theta_\ell$ is near zero $\theta = 0, ~2 \pi$.  Thus for $\theta_\ell$ away from zero, the integral will be of lower order than for $\theta_\ell$ near $\theta = 0, 2\pi$.  The above calculations capture this fact.
\end{remark}

We use the notation $\mathcal R_nf'$ to refer to the derivative of $\mathcal R_nf$.  Examine
\begin{align*}
\int_{\mathbb R} |\mathcal R_nf'(x)|^2 dx = \int_{0}^{2\pi} |\mathcal I_nF'(\theta)|^2 |T'(\theta)|^{-1} d\theta.
\end{align*} 
As before, the factor $|T'(\theta)|^{-1}$ vanishes at $\theta = 0, 2\pi$ which makes the bounding of the operator easier.  Proceeding,
\begin{align*}
\|\mathcal R_nf'\|_{L^2(\mathbb R)} \leq \|F\|_{C_p^0(\mathbb R)}\|D_n' |T'|^{-1}\|_{L^2(\mathbb T)} \leq C \|F\|_{C_p^0(\mathbb R)}\|D_n'\|_{L^2(\mathbb T)}
\end{align*}
Combining these estimates with Theorem~\ref{Theorem:Galeev}, we find
\begin{align*}
\|\mathcal R_nf\|_{H^1(\mathbb R)} \leq C n^{3/2}\|f\|_{C_0^0(\mathbb R)}\leq C_{1} n^{3/2}\|f\|_{H^1(\mathbb R)}.
\end{align*}

\begin{remark}
Because $D'_n(\theta-\theta_\ell)$ is largest near $\theta_\ell$, one might expect the vanishing of $|T'(\theta)|^{-1}$ at $\theta = 0,~2\pi$ to reduce the magnitude of the integral for $\theta_\ell$ near $0,~2\pi$.  While this does indeed happen, once the sum over $\ell$ is performed, the result is still $\bigo(n^{3/2})$.
\end{remark}

\begin{remark}
In proving a bound on the $H^1(\mathbb R) \goto H^1(\mathbb R)$ operator norm we passed through $C_0^0(\mathbb R)$.  Presumably, a tighter bound can be found by using further structure of the rational approximation of an $H^1(\mathbb R)$ function. One will no longer be able to make use of the $H^s(\mathbb T)$ theory and therefore refining this estimate is beyond the scope of the current paper.
\end{remark}

\section{Oscillatory Cauchy integrals}\label{Section:OscCauchy}

At this point, we have developed the algorithm and theory for a method that provides a rapidly convergent rational approximation of $f : \mathbb R \goto \mathbb C$ provided that $f(T(\theta))$ is sufficiently smooth.  This approximation can, depending on the amount of smoothness, converge in a whole host of Sobolev spaces. We now review a formula from \cite{Trogdon2013} for the computation of Cauchy integrals of the form
\begin{align}\label{Cauchy-int}
\frac{1}{2\pi i} \int_{\mathbb R} e^{-i k x}\left[\left(\frac{x-i\beta}{x+i\beta}\right)^j -1 \right] \frac{\D x}{x-z}, ~~ z \in \mathbb C \setminus \mathbb R.
\end{align}
In other words, we compute the Cauchy integral of the oscillatory rational basis $\{R_{j,k}\}_{j=-\infty}^\infty$, $k \in \mathbb R$ where $R_{j,k}(x) = e^{-ikx}R_j(x)$. The notation here differs from that in \cite{Trogdon2013} by a sign change of $k$. We use the notation $\mathcal C_{\mathbb R} R_{j,k}(z)$ to denote \eqref{Cauchy-int}.  Furthermore,
\begin{align*}
\mathcal C_{\mathbb R}^\pm R_{j,k}(x) = \lim_{\epsilon \goto 0^+}  \mathcal C_{\mathbb R} R_{j,k}(x \pm i \epsilon) 
\end{align*}
are used to denote the boundary values.

\begin{remark}\label{Remark:mult}
This basis is closed under pointwise multiplication.  A straightforward calculation shows the simple relation for $k_1,~k_2 \in \mathbb R$
\begin{align*}
R_{j,k_1}(z)R_{\ell,k_2}(z) = R_{\ell+j,k_1+k_2}(z) -R_{j,k_1+k_2}(z)-R_{\ell,k_1 + k_2}(z).
\end{align*}
\end{remark}

 An important aspect of the formula we review is that it expresses  $\mathcal C_{\mathbb R} R_{j,k}(z)$ in terms of $R_{j,k}(x)$ and $R_{j,0}(x)$ so that it is useful for the approximation of operator equations \cite{Trogdon2013}.  Define for $j > 0$, $n > 0$,
\begin{align*}
\gamma_{j,n}(k) = - \frac{j}{n} e^{- |k| \beta} \left(\begin{array}{cc} j-1 \\ n \end{array}\right) \phantom{.}_1F_1(n-j,1+n,2 |k| \beta),
\end{align*}
where $_1F_1$ is Krummer's confluent hypergeometric function
\begin{align*}
_1F_1(a,b,z) = \sum_{\ell =0}^\infty \frac{\Gamma(a +\ell)}{\Gamma(a)} \frac{\Gamma(b)}{\Gamma(b+\ell)} \frac{z^\ell}{\ell!},
\end{align*}
and $\Gamma$ denotes the Gamma function \cite{DLMF}.   Further, define
\begin{align*}
\eta_{j,n}(k) = \sum_{\ell=n}^{j} (-1)^{n+\ell} \left(\begin{array}{cc} \ell \\ n \end{array}\right) \gamma_{\ell,n}(k).
\end{align*}
The following theorem is proved by pure residue calculations.
\begin{theorem}[\cite{Trogdon2013}]\label{Theorem:CauchyAction}
If $k j < 0$ then
\begin{align*}
\mathcal C^+_{\mathbb R} R_{j,k}(x) &= \begin{choices} R_{j,k}(x), \when j > 0,\\ 0, \when j < 0, \end{choices}\\
\mathcal C^-_{\mathbb R} R_{j,k}(x) &= \begin{choices} 0, \when j > 0, \\ -R_{j,k}(x), \when j < 0.\end{choices}
\end{align*}
If $k j > 0$ then
\begin{align*}
\mathcal C^+_{\mathbb R} R_{j,k}(x) & = \begin{choices} - \displaystyle \sum_{n=1}^{j} \eta_{j,n}(k) R_{n,0}(x), \when j > 0,\\ 
R_{j,k}(x) +  \displaystyle\sum_{n=1}^{-j} \eta_{j,n}(k) R_{-n,0}(x), \when j < 0, \end{choices}\\
\mathcal C^-_{\mathbb R} R_{j,k}(x) & = \begin{choices} -R_{j,k}(x)-  \displaystyle\sum_{n=1}^{j} \eta_{j,n}(k) R_{n,0}(x), \when j > 0,\\ 
  \displaystyle\sum_{n=1}^{-j} \eta_{j,n}(k) R_{-n,0}(x), \when j < 0. \end{choices}
\end{align*}
\end{theorem}

Define the oscillatory Cauchy operator
\begin{align*}
\mathcal C_{\mathbb R,k} f(z) = \frac{1}{2\pi i} \int_{\mathbb R} \frac{e^{-ikx}f(x)}{x-z} \D x.
\end{align*}
Next, define the approximation of the Cauchy integral of $e^{-ikx}f(x)$ for $f \in H^1(\mathbb R)$
\begin{align*}
\mathcal C_{\mathbb R,k,n}f(z) = \mathcal C_{\mathbb R,k} \mathcal R_{n}f(z),
\end{align*}
and its boundary values
\begin{align*}
\mathcal C^\pm_{\mathbb R,k,n}f(x) = \mathcal C^\pm_{\mathbb R,k} \mathcal R_{n}f(x).
\end{align*}

\subsection{Accuracy}

We now address the accuracy of the operator $\mathcal C_{\mathbb R,n,k}f(z)$ both on the real axis and in the complex plane.  To do this, we require some well-known results concerning the Cauchy integral operator (see, for example, \cite{DeiftOrthogonalPolynomials,wavelets,TrogdonThesis}).

\begin{theorem}\label{Theorem:CauchyProp}
Assume $f \in H^s(\mathbb R)$, $s \in \mathbb N$ and $\delta > 0$ then
\begin{itemize}
\item  $\mathcal C_{\mathbb R} f( \cdot \pm i \delta) \in H^s(\mathbb R)$, $\|\mathcal C_{\mathbb R} f(\cdot \pm i \delta) \|_{H^s(\mathbb R)} \leq \|\mathcal C_{\mathbb R}^\pm f\|_{H^s(\mathbb R)} \leq \|f\|_{H^s(\mathbb R)}$ and
\item for $\Omega_\delta = \{\pm (z + i \delta) ~:~ \imag z > 0 \}$, $\sup_{z \in \Omega_{\delta}} |\mathcal C_{\mathbb R} f^{(j)}(z)| \leq C_{j,\delta} \|f\|_{L^2(\mathbb R)}$.
\end{itemize}
\end{theorem}

\begin{remark}
The first statement of this theorem follows directly from the fact that the Fourier symbol for the Cauchy integral operators is bounded by unity and therefore it does not destroy $H^s(\mathbb R)$ smoothness.
\end{remark}

A straightforward calculation produces
\begin{align*}
\|\mathcal C_{\mathbb R,k} f\|_{H^s(\mathbb R)} \leq C_s\sum_{j=0}^s |k|^{s-j} \|f\|_{H^{j}(\mathbb R)}.
\end{align*}
Therefore, we must be aware that $H^s(\mathbb R)$ errors made in the approximation of $f$ may be amplified as $|k|$ increases.  We concentrate on $L^2$, $H^1$ and uniform convergence in what follows.

Two inequalities easily follow from Theorems~\ref{Theorem:Convergence} and \ref{Theorem:CauchyProp}, uniform for $\delta \geq 0$,
\begin{align*}
\|(\mathcal C_{\mathbb R,k} - \mathcal C_{\mathbb R,k,n})f\|_{H^1(\mathbb R+ i \delta)} &\leq (1+|k|)\|\mathcal R_nf-f\|_{H^1(\mathbb R)}\\
&\leq C_{\epsilon,s}(1+|k|) n^{1+ \epsilon -s} \|F\|_{H^s(\mathbb T)}\\
\|(\mathcal C_{\mathbb R,k} - \mathcal C_{\mathbb R,k,n})f\|_{L^2(\mathbb R+ i \delta)} &\leq \|\mathcal R_nf-f\|_{L^2(\mathbb R)}\\
&\leq C_{\epsilon,s}n^{1+ \epsilon -s} \|F\|_{H^s(\mathbb T)},
\end{align*}
for $\delta > 0$ and $F(\theta) = f(T(\theta))$, as before.  Then, for $\delta  > 0$ by Sobolev embedding and Theorem~\ref{Theorem:CauchyProp}
\begin{align*}
\sup_{x \in \mathbb R} |\mathcal C_{\mathbb R,k}f(x \pm i \delta)  &- \mathcal C_{\mathbb R, k,n}f(x \pm i \delta)| \\
&\leq \|(\mathcal C_{\mathbb R,k} - \mathcal C_{\mathbb R,k,n})f\|_{H^1(\mathbb R)} \leq C_{\epsilon,s}(1+|k|) n^{1+ \epsilon -s} \|F\|_{H^s(\mathbb T)},
\end{align*}
and for fixed $k$, we realize uniform convergence on all of $\mathbb C$.

\begin{remark}
This theoretical result can be a bit misleading.  Consider computing $\mathcal C^+_{\mathbb R} R_{j,k}$ when $j < 0$ and $k > 0$.  From Theorem~\ref{Theorem:CauchyAction}
\begin{align}\label{cp-odd}
\mathcal C^+_{\mathbb R} R_{j,k}(z)  = R_{j,k}(z) + \sum_{n=1}^{-j} \eta_{j,n}(k) R_{-n,0}(z).
\end{align}
Necessarily, the operator $\mathcal C^+_{\mathbb R}$ produces a function that is analytic in the upper-half of the complex plane.  Each term in \eqref{cp-odd} has a pole at $z = i$!  Very specific cancellation occurs to ensure that this function is analytic at $z = i$.  Therefore, the evaluation of this formula for large $j$ is not stable.  There is a similar situation for $\mathcal C_{\mathbb R}^-R_{j,k}$ when $j > 0$ and $ k < 0$.
\end{remark}


\section{Fourier integrals}\label{Section:Fourier}

The same methods for computing oscillatory Cauchy integrals apply to the computation of Fourier transforms of functions that are well-approximated in the basis $\{R_j\}_{j=-\infty}^\infty$.  The relevant expression from \cite{Trogdon2013} is
\begin{align*}
\mathcal F{R}_j(k) = \dashint_{\mathbb R} R_{j,k}(z) dz = \omega_j(k) = \begin{choices}
0, \when \sign(j) = -\sign(k),\\
- 2 \pi |j| \beta, \when k = 0,\\
\displaystyle -4 \pi e^{- |k| \beta} \beta L_{|j|-1}^{(2)}(2 |k|\beta), \otherwise,\end{choices}
\end{align*}
where $L^{(\alpha)}_n(x)$ is the generalized Leguerre polynomial of order $n$ \cite{DLMF}.  A similar expression was discussed in \cite{Henrici1979} for a slightly different rational basis in the context of the Laplace transform (see also \cite{Weber1980}). The $L^2(\mathbb R)$ convergence of $\mathcal F \mathcal R_n f$ is easily analyzed with the help of Definition~\ref{Def:SobolevLine} and Theorem~\ref{Theorem:Convergence}.  For $F = f \circ T \in H^s(\mathbb T)$ and $t \leq s$
\begin{align*}
\|\mathcal F(\mathcal R_n f -f)(1+|\cdot|)^{t}\|_{L^2(\mathbb R)} = \|\mathcal R_n f -f\|_{H^t(\mathbb R)} \leq [C_{\epsilon,s} n^{1+\epsilon-s} + C_{s,t} n^{t-s}] \|F\|_{H^s(\mathbb T)}.
\end{align*}

We remark that because generically, $\mathcal R_n f(x) = \bigo(|x|^{-1})$ as $x \goto \infty$, $\mathcal F\mathcal R_nf(k)$ will have a discontinuity at the origin.  This approximation will not converge in any Sobolev space.  It will converge uniformly as we now discuss.  Assuming $F$ is continuously differentiable, let $s_n$ be as in \eqref{fixed-pv} and define
\begin{align*}
S_n(f) &= \left| \dashint_{\mathbb R} e^{-ikx}(\mathcal R_nf(x) - f(x)) dx \right| \\
&\leq \left|\int_{\mathbb R} e^{-ikx}(\mathcal R_nf(x) + i s_n R_1(x) - f(x)) dx  \right| + \left|s_n \dashint e^{-ikx}R_1(x) dx \right|,\\
& \leq \int_{\mathbb R} |\mathcal R_nf(x) + i s_n R_1(x) - f(x))| dx + 4 \pi |s_n| e^{- |k| \beta} \beta.
\end{align*}
Following the arguments in Section~\ref{Section:ConvIntegral} the first term is bounded by $C_{\epsilon,s} n^{3/2+\epsilon -s } \|F\|_{H^s(\mathbb T)}$ and $|s_n| \leq C_s n^{3/2-s}\|F\|_{H^s(\mathbb T)}$ so that
\begin{align*}
|\mathcal F \mathcal R_n f(k) - \mathcal F f(k)| = S_n(f) \leq C_{\epsilon,s} n^{3/2+\epsilon-s} \|F\|_{H^s(\mathbb T)},~~ s > 3/2 + \epsilon,
\end{align*}
proving uniform convergence of the Fourier transforms.  More is true. Performing integration by parts $j$ times on $S_n(f)$ we find for $|k| > 0$ and $j > 0$
\begin{align*}
S_n(f) \leq |k|^{-j} \| \mathcal D^j [\mathcal R_n f - f]\|_{L^1(\mathbb R)} \leq C_{j,s} |k|^{-j} n^{j-s} \|F\|_{H^s(\mathbb T)},
\end{align*}
from Theorem~\ref{Theorem:Convergence}.  Hence we realize an asymptotic approximation of Fourier integrals, for fixed $n$.

\subsection{An oscillatory quadrature formula}

From above, the coefficients for the rational approximation in the basis $\{R_{j}(x)\}_{j\in \mathbb Z}$ are given by
\begin{align*}
\tilde F_j = \frac{1}{n} \sum_{\ell = 0}^{n-1} e^{ij\theta_\ell} f(M_{\beta}^{-1}(e^{i\theta_\ell})).
\end{align*}
The approximation of the Fourier transform is given by
\begin{align*}
\mathcal Ff(k) \approx \sum_{j=-n_-}^{n^+} \omega_j(k) \frac{1}{n} \sum_{\ell = 0}^{n-1} e^{ij\theta_\ell} f(M_{\beta}^{-1}(e^{i\theta_\ell})) = \sum_{\ell=0}^{n-1} f(M_{\beta}^{-1}(e^{i\theta_\ell})) \left( \sum_{j=-n_-}^{n_+} \frac{\omega_j(k)}{n} e^{ij\theta_\ell}  \right)
\end{align*}
Therefore we have the interpretation of this formula with quadrature nodes $\{M_{\beta}^{-1}(e^{i\theta_\ell})\}_{\ell=0}^{n-1}$ and weights
\begin{align*}
w_\ell = \sum_{j=-n_-}^{n_+} \frac{\omega_j(k)}{n} e^{ij\theta_\ell}.
\end{align*}
Despite this interpretation, it is more efficient to compute the coefficients with the fast Fourier transform and perform a sum to compute the Fourier integral.  This is true, unless, this sum can be written in closed form.

\section{Differentiation}

The differentation operator $\mathcal D$ also maps the basis $\{R_{j,k}\}_{j\in \mathbb Z, ~k \in \mathbb R}$ to itself in a convenient way.  Based on the theory we have developed, we know precise conditions to impose on $f$ so that $\mathcal D \mathcal R_{n} f$ converges to $\mathcal Df$ uniformly, in $L^2(\mathbb R)$ and in $L^1(\mathbb R)$.  For $j >0$ we consider
\begin{align}\label{diff-Rj}
R_j'(z) = \frac{d}{dz} \left[ \left( \frac{z - i \beta}{z+i\beta} \right)^j -1 \right] = j\left( \frac{z - i \beta}{z+i\beta} \right)^{j-1} \frac{2 i \beta}{(z+i\beta)^2}.
\end{align}
A simple computation shows that
\begin{align*}
\frac{2 i \beta}{(z+i\beta)^2} = i \frac{1}{\beta} R_1(z) - i \frac{1}{2\beta} R_2(z).
\end{align*}
Furthermore, following Remark~\ref{Remark:mult} we find
\begin{align*}
R_j'(z) &= j R_{j-1}(z) \left( i \frac{1}{\beta} R_1(z) - i\frac{1}{2\beta} R_2(z) \right) + i \frac{j}{\beta} R_1(z) - i\frac{j}{2\beta} R_2(z),\\
&= i \frac{j}{\beta} R_j(z) - i \frac{j}{\beta} R_{j-1}(z) - i\frac{j}{\beta} R_1(z) - i\frac{j}{2\beta} R_{j+1}(z) \\
&+ j\frac{i}{2\beta} R_{j-1} + j\frac{i}{2\beta} R_{2}(z) + i \frac{j}{\beta} R_1(z) - i\frac{j}{2\beta} R_2(z)\\
&= - i\frac{j}{2\beta} R_{j+1}(z) + i\frac{j}{\beta} R_j(z)  - i\frac{j}{2\beta}R_{j-1}(z).
\end{align*}
To obtain a formula for $j < 0$, note that $R_j'(z)$ is just the complex conjugate of \eqref{diff-Rj} for $z \in \mathbb R$.  Thus ($j < 0$)
\begin{align*}
R_j'(z) = \overline{R_{-j}'(z)} =  i\frac{j}{2\beta} R_{j-1}(z) - i\frac{j}{\beta} R_{j}(z)  + i\frac{j}{2\beta}R_{j+1}(z),
\end{align*}
where we used that $\overline{R_j(z)} = R_{-j}(z)$ for $z \in \mathbb R$.  This is summarized in the following proposition.

\begin{proposition}
Suppose $f \in C_0^0(\mathbb R)$, and
\begin{align*}
\mathcal R_n f(x) = \sum_{j=-n_-}^{n_-} \alpha_j R_j(x)
\end{align*}
then
\begin{align*}\mathcal D \mathcal R_n f(x) = \sum_{j=-n_-}^{n_-} \frac{i}{\beta} \left( - |j-1|\alpha_{j-1} + |j| \alpha_j  -|j+1|\alpha_{j+1}\right) R_j(x),
\end{align*}
\emph{i.e.}
\begin{align*}
R_j'(z) = \sign(j) \left( -i\frac{j}{2\beta} R_{j-1}(z) + i\frac{j}{\beta} R_{j}(z)  - i\frac{j}{2\beta}R_{j+1}(z) \right).
\end{align*}
Note that all sums can be taken to avoid the index $j =0$ ($R_0(z) = 0$) and we use the convention that $\alpha_j = 0$ for $|j|>n$.
\end{proposition}

\section*{Acknowledgments}

The author would like to thank Sheehan Olver for discussions that led to this work.  The author acknowledges the National Science Foundation for its generous support through grant NSF-DMS-130318.  Any opinions, findings, and conclusions or recommendations expressed in this material are those of the author and do not necessarily reflect the views of the funding sources.

\appendix

\section{A Numerical Example}

\begin{figure}[tbp]
\centering
\includegraphics[width=.9\linewidth]{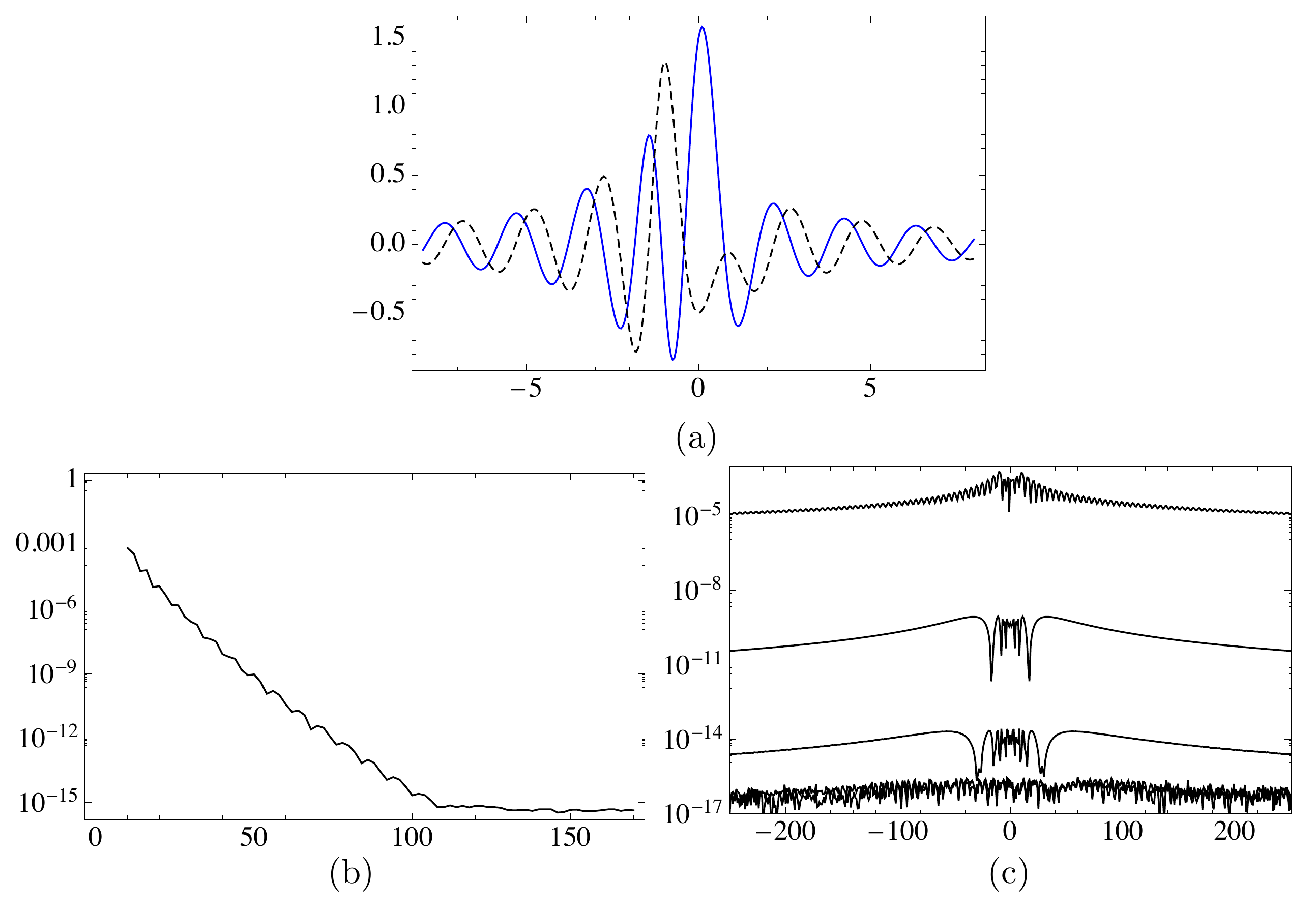}
\caption{\label{Fig:Rnf} (a) A plot of the function $f$ with $k_1 = 2,~k_2 = -3$ (solid: real part, dashed: imaginary part).  (b) A convergence plot of an estimate of $\sup_{|x|\leq 60} |e^{-ik_1x}\mathcal R_nf_1(x)+ e^{-ik_2x}\mathcal R_nf_2(x) - f(x)|$ plotted on a log scale versus $n$.  This plot clearly shows super-algebraic (spectral) convergence.  (c) A plot of $|e^{-ik_1x}\mathcal R_nf_1(x)+ e^{-ik_2x}\mathcal R_nf_2(x) - f(x)|$ versus $x$ for $n = 10, ~50, ~90, ~130$.}
\end{figure}

\begin{figure}[tbp]
\centering
\includegraphics[width=.9\linewidth]{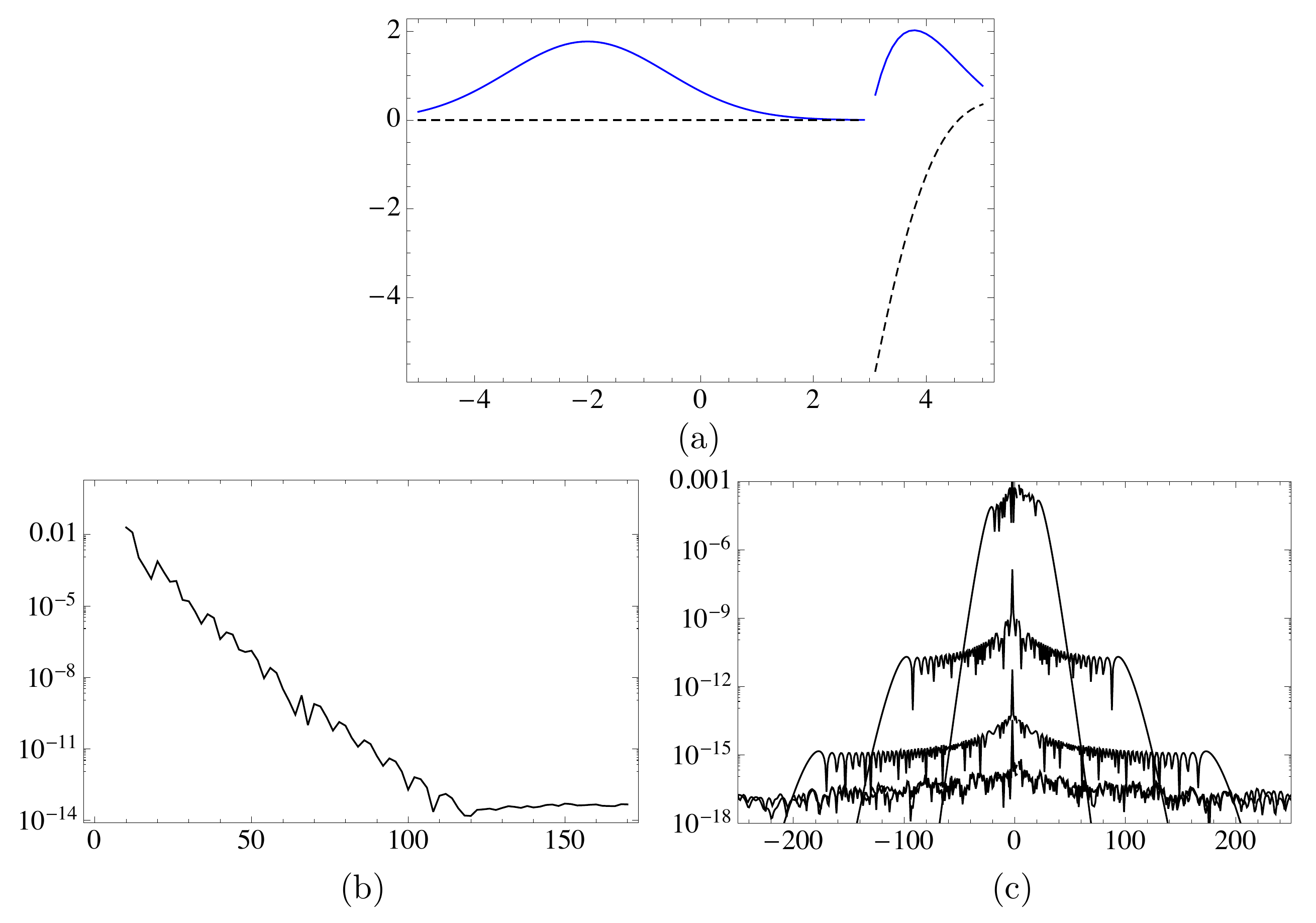}
\caption{\label{Fig:FRnf}(a) A plot of $\mathcal Ff$, the Fourier transform of $f$, with $k_1 = 2,~k_2 = -3$ (solid: real part, dashed: imaginary part).  Note that this function has a discontinuity at $k = 3$.  (b) A convergence plot of an estimate of $\sup_{|x|\leq 60} |\mathcal F[e^{-ik_1x}\mathcal R_nf_1(x)+ e^{-ik_2x}\mathcal R_nf_2(x)] - \mathcal Ff(x)|$ plotted on a log scale versus $n$.  This plot clearly shows super-algebraic (spectral) convergence.  (c) A plot of $|\mathcal F[e^{-ik_1x}\mathcal R_nf_1(x)+ e^{-ik_2x}\mathcal R_nf_2(x)] - \mathcal Ff(x)|$ versus $x$ for $n = 10, ~50, ~90, ~130$.  The exponential decay of the Fourier transforms of $R_j$ is evident.  In this way, the method is accurate asymptotically.}
\end{figure}

\begin{figure}[tbp]
\centering
\includegraphics[width=.9\linewidth]{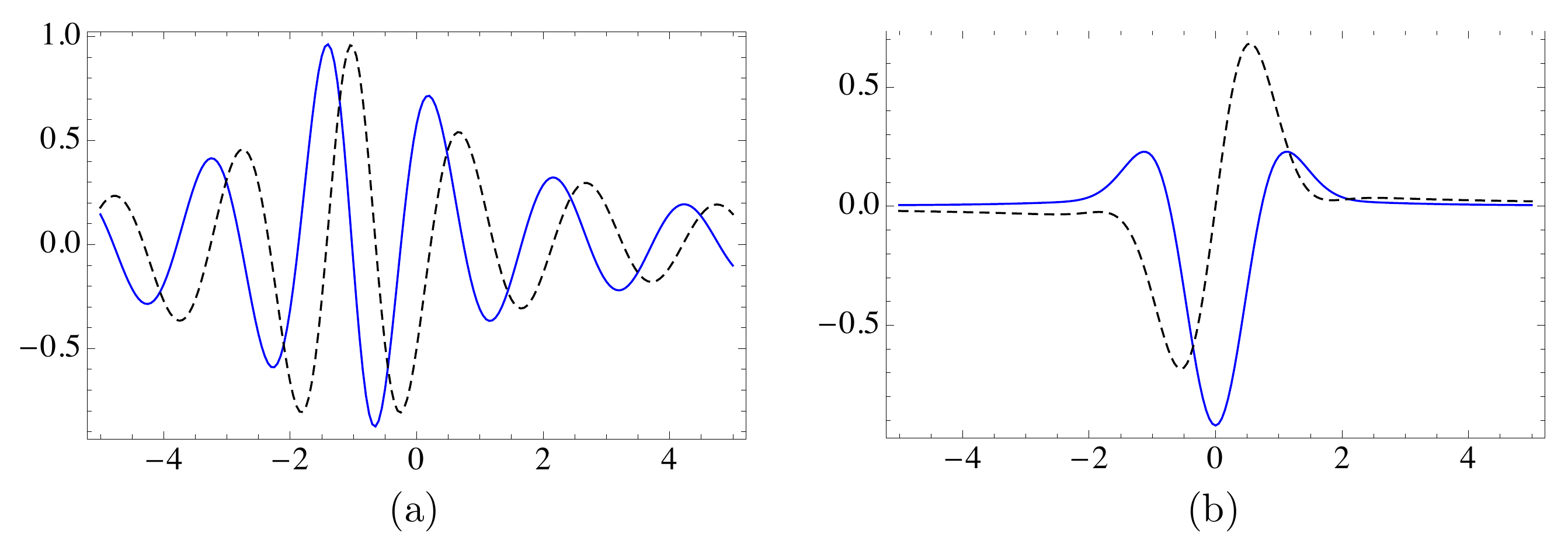}
\caption{\label{Fig:CpmRnf} (a) The Cauchy transform $\mathcal C^+_{\mathbb R}f(x)$ with $k_1 =2,~k_2 = -3$ plotted on the real axis (solid: real part, dashed: imaginary part). (b) The Cauchy transform $\mathcal C^-_{\mathbb R}f(x)$ plotted on the real axis (solid: real part, dashed: imaginary part).  Note that due to the fact that $k_1 > 0$, the $\mathcal C^-_{\mathbb R}$ operator essentially isolates the Gaussian term in $f$.  Also, because $k_2 < 0$, the rational term in $f$ is isolated by $\mathcal C_{\mathbb R}^+$.}
\end{figure}

\begin{figure}[tbp]
\centering
\includegraphics[width=.9\linewidth]{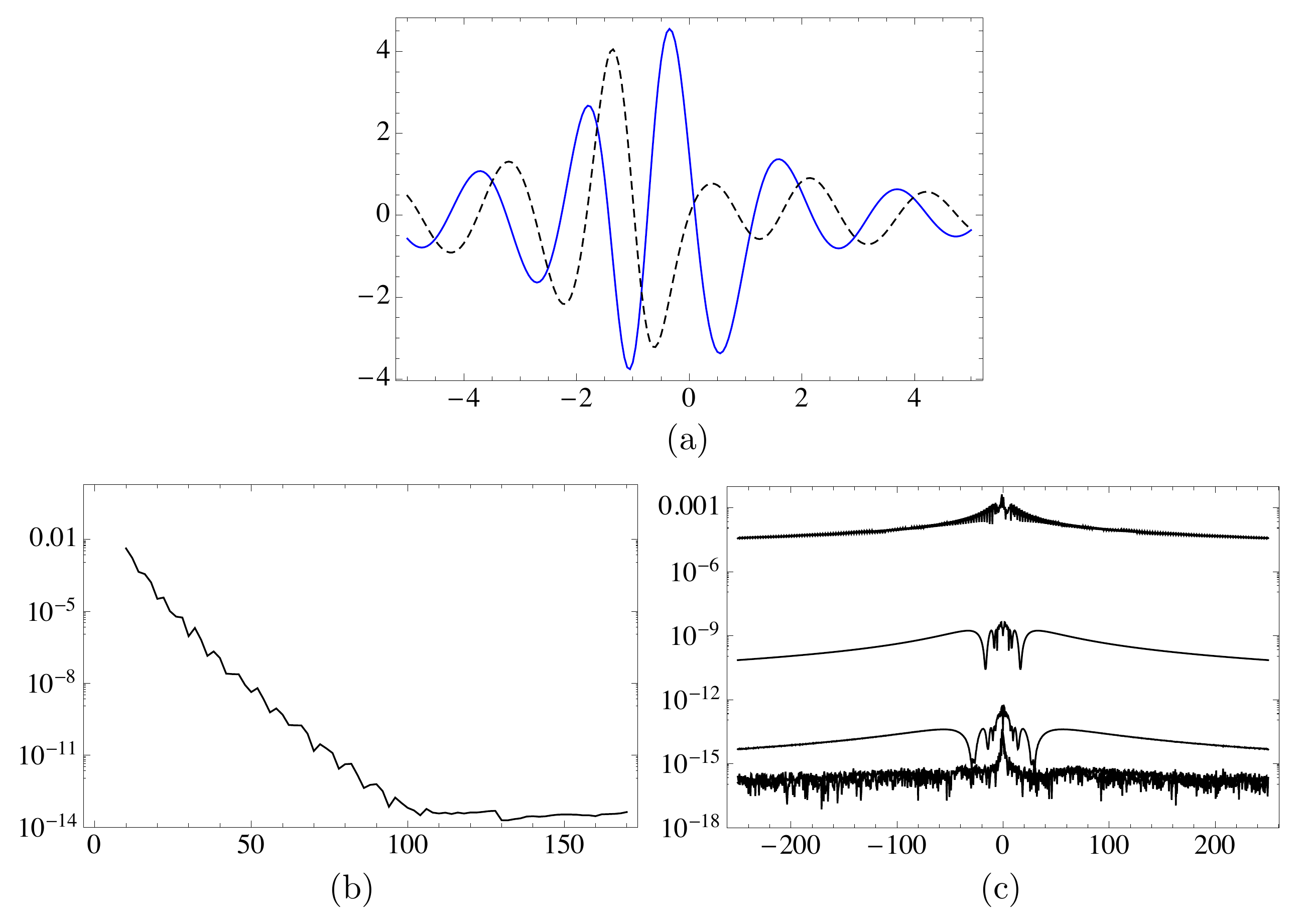}
\caption{\label{Fig:DRnf}(a) A plot of $\mathcal Df$, the derivative of $f$, with $k_1 = 2,~k_2 = -3$ (solid: real part, dashed: imaginary part).  (b) A convergence plot of an estimate of $\sup_{|x|\leq 60} |\mathcal D[e^{-ik_1x}\mathcal R_nf_1(x)+ e^{-ik_2x}\mathcal R_nf_2(x)] - \mathcal Df(x)|$ plotted on a log scale versus $n$.  Just like the previous plot, this clearly shows super-algebraic (spectral) convergence.  (c) A plot of $|\mathcal D[e^{-ik_1x}\mathcal R_nf_1(x)+ e^{-ik_2x}\mathcal R_nf_2(x)] - \mathcal Df(x)|$ versus $x$ for $n = 10, ~50, ~90, ~130$.  Note that while absolute error is lost (when comparing with Figure~\ref{Fig:Rnf}) near the origin, litte accuracy is lost, if any, in the tails.  This is related to the phenomenon that is described in Theorem~\ref{Theorem:Convergence}: approximation of the function in $H^1(\mathbb R)$ is almost as ``easy'' as approximation in $L^2(\mathbb R)$.}
\end{figure}

To demonstrate the method we consider the function
\begin{align*}
f(x) = e^{-ik_1x -x^2} + \frac{e^{-ik_2x}}{x + i+1}.
\end{align*}
This function can clearly be written in the form
\begin{align*}
f(x) = e^{-ik_1x}f_1(x) + e^{-ik_2x} f_2(x),
\end{align*}
where $f_1$ and $f_2$ have rapidly converging interpolants (see Propositions~\ref{Prop:Decay} and \ref{Prop:Analytic}).  We demonstrate the rational approximation of $f$ in Figure~\ref{Fig:Rnf}.  The approximation of the Fourier transform of $f$, which contains a discontinuity, is demonstrated in Figure~\ref{Fig:FRnf}.  The Cauchy operators applied to $f$ are shown in Figure~\ref{Fig:CpmRnf}.  Finally, the convergence of the approximation of the derivative of $f$ is shown in Figure~\ref{Fig:DRnf}.

\section{Sufficient Conditions for Spectral Convergence}

To demonstrate spectral convergence it is sufficient to have $F \in H^s(\mathbb T)$ for every $s > 0$.  We provide two propositions that provide sufficient conditions for this.
\begin{proposition}\label{Prop:Decay}
Suppose $f^{(\ell)}(x)$ exists and is continuous for every $\ell > 0$.  Further, assume $\sup_{x \in \mathbb R}(1+|x|)^j|f^{(\ell)}(x)|<\infty$ for every $j,\ell \in \mathbb N$. Then $F \in H^s(\mathbb T)$ for every $s > 0$.
\begin{proof}
Every derivative of $f(x)$ decays faster then any power of $x$.  It is clear that all derivatives of $F(\theta) = f(T(\theta))$ are continuous at $\theta = 0, ~2\pi$.
\end{proof}
\end{proposition}
From this proposition it is clear that the method will converge spectrally for, say, $f(x) = e^{-x^2}$.

\begin{proposition}\label{Prop:Analytic}
Suppose $f$ can be expressed in the form
\begin{align*}
f(x) = \frac{1}{2 \pi i} \int_{\Gamma} \frac{h(s)}{s-x} ds
\end{align*}
for some oriented contour $\Gamma \subset \mathbb C$ and $|x|^jf(x) \in L^1(\Gamma)$ for every $j \geq 0$.  Assume further that $\min_{x \in \mathbb R, s \in \Gamma} |x-s| \geq \delta > 0$ then $F \in H^s(\mathbb T)$ for every $s$.
\begin{proof}
We claim that it is sufficient to prove that there exists an asymptotic series
\begin{align*}
f(x) \sim \sum_{j=1}^\infty c_j x^{-j}, \quad x \goto \pm \infty,
\end{align*}
that can be differentiated (see \cite{SOHilbertTransform}).  A short computation shows that this is sufficient for $F(\theta)$ to be smooth in a neighborhood of $\theta = 0,~2\pi$.  Furthermore, it is clear that $f$ is smooth on $\mathbb R$.  Thus the differentiable asymptotic series is sufficient for $F \in H^s(\mathbb T)$ for all $s$.

Our guess for the asymptotic series is obtained from a Neumann series expansion for $1/(s-x)$ for $s$ fixed and $|z|$ large:
\begin{align*}
\frac{1}{s-x} = - \frac{1}{x} \sum_{j=0}^{n-1} \left( \frac{s}{x} \right)^j - \frac{(s/x)^n}{s-x}.
\end{align*}
Therefore we choose
\begin{align*}
c_j = - \frac{1}{2\pi i} \int_{\Gamma} s^{j-1} f(s) ds,
\end{align*}
and define $P_n(s,x) = (s/x)^n/(x-s)$. We find for $x \in \mathbb R$, $|x| \geq 1$,
\begin{align*}
\sup_{s \in \Gamma} |\partial_x^j P_n(s,x)| \leq C_{\delta,j,n}\left| \frac{s}{x} \right|^n.
\end{align*}
This estimate dictates the order of the error term and justifies the interchange of differentation and integration.  Hence
\begin{align*}
f^{(\ell)}(x) - \sum_{j=1}^{n-1} (-1)^\ell c_j& \frac{(j+\ell-1)!}{(j-1)!} x^{-j-\ell}  \\
&= \sum_{j=n}^{n+\ell-1} c_j \frac{(j+\ell-1)!}{(j-1)!} x^{-j-\ell} + \frac{1}{2 \pi i} \int_{\Gamma} \partial_x^\ell P_{n+\ell}(s,x)f(s) ds = \bigo(x^{-n-\ell}).
\end{align*}
This shows the existence and differentiability of the asymptotic series.
\end{proof}
\end{proposition}
Any rational function that decays at infinity can be expressed in terms of a contour integral in the finite complex plane.  This proposition, in particular, demonstrates spectral convergence for such functions.

\section{$L^p$ Estimates of the Dirichlet Kernel}\label{Appendix:LpEstimates}

Our first step is to rewrite the Dirichlet kernel
\begin{align*}
D_n(\theta) = \sum_{k=-n_-}^{n_+} e^{ik\theta} = \sum_{k=-n_-}^{n_-} e^{ik\theta} + \sigma e^{ikn_+},
\end{align*}
where $\sigma = 1$ if $n_+ > n_1$ ($n$ is even) and $\sigma = 0$ otherwise. It is well known that this may be expressed as
\begin{align*}
D_n(\theta) = \frac{\sin( (n_-+1/2) \theta)}{2 \sin \theta/2} + \sigma e^{ikn_+}.
\end{align*}

It is clear that the $L^p(\mathbb T)$ norm of $e^{ikn_+}$ is uniformly bounded by a constant.  Therefore, the large $n$ behavior is dictated by the ratio of sines.  Consider the integral for $\pi \geq \epsilon/n_- \geq 0$, and using periodicity
\begin{align*}
I(n,p) &= \int_{\epsilon/n_-}^{2\pi-\epsilon/n_-} \left| \frac{\sin( (n_-+1/2) \theta)}{2 \sin \theta/2} \right|^p d \theta = \left( \int_{-\pi}^{-\epsilon/n_-} + \int_{\epsilon/n_-}^{\pi}\right)\left| \frac{\sin( (n_-+1/2) \theta)}{2 \sin \theta/2} \right|^p d \theta.
\end{align*}
Let $y = \theta n_-$ so that
\begin{align*}
I(n,p) n_-^{1-p} =  \left( \int_{-\pi n_-}^{-\epsilon} + \int_{\epsilon}^{\pi n_-}\right)\left| \frac{\sin( (1+1/(2n_-) y) \theta)}{2 n_- \sin (y/(2n_-))} \right|^p d y.
\end{align*}
We bound the ratio of sines.  First,
\begin{align*}
|\sin(y+y/(2n_-))| &= | \sin y \cos (y/(2n_-))  + \cos y \sin(y/(2 n_-))| \\
&\leq |\sin y | + |\cos y| \min\{y/(2n_-),1\}.
\end{align*}
Next, because $|\sin y | \leq |y|$ and the Taylor series for $\sin y$ is an alternating series
\begin{align*}
|y| - |\sin y| \leq |y|^3/6 ~~\Rightarrow~~ |y|(1 - |y|^2/6) \leq |\sin y|.
\end{align*}
This bound is useful provided the left-hand side stays positive: $|y| < \sqrt{6}$.  We find
\begin{align*}
\frac{1}{|2n_- \sin(y/(2n_-))|} \leq \frac{1}{|y|} \max_{y \in [-\pi n_-, \pi n_-]} (1 - |y|^2/(24n_-^2))^{-1} \leq  \frac{2}{|y|}
\end{align*}
Our inequalities demonstrate that
\begin{align*}
I(n,p) n_-^{1-p} &\leq 4\int_{\epsilon}^\infty \left|\frac{|\sin y | + |\cos y| \min\{y/(2n_-),1\}}{y}\right|^p dy \leq C_p^p (1+\epsilon)^{1-p},\\
\|D_n\|_{L^p(\epsilon/n_-,2 \pi - \epsilon/n_-)} &\leq C_p \left(\frac{n_-}{1+\epsilon}\right)^{1-1/p}.
\end{align*}
Also, setting $\epsilon = 0$ results in a special case of Theorem~\ref{Theorem:Galeev}.

We now turn to the first derivative and perform calculations in the case that $n$ is odd so that $\sigma = 0$.  If $n$ were even then the $L^p(\mathbb T)$ norm of the derivative of $e^{ik n_-}$ is $\bigo (n)$.  If we restrict to $\pi \geq \epsilon/n_- \geq 0$ then the bound we obtain on rest of $D_n$ is of equal or larger order, justifying setting $\sigma = 0$.  Consider
\begin{align*}
D_n'(\theta) = \frac{2 (n_-+1/2)\cos ((n_-+1/2) \theta) \sin \theta/2 - \cos\theta/2 \sin ((n_-+1/2) \theta)}{4 \sin^2 \theta/2}.
\end{align*}
Again, we invoke the change of variables $y = \theta n_-$.  The same procedure as above, in principle works.  We have already derived an appropriate bound on the denominator after we divide by $n_-^2$.  To bound the numerator, we look for a bound, independent of $n$ that vanishes to second order at $y = 0$ to cancel the singularity from the denominator.  We show that there exists $C > 0$ so that
\begin{align}\label{Dprime-est}
\frac{1}{n_-^2} |D_n'(y/n_-)| \leq C(1 + |y|)^{-1}, ~~ y \in [-\pi n_-, \pi n_-].
\end{align}
From above we have
\begin{align}\label{denom-bound}
|4 n_-^2 \sin^2 (y/(2n_-)) |^{-1} \leq 4 |y|^{-2}
\end{align}
so that we just need to estimate the numerator.  Straightforward trigonometric rearrangements show that
\begin{align}
2 (n_-+1/2)\cos (y + y/(2n_-) ) \sin (y/(2n_-)) &- \cos(y/(2n_-)) \sin (y + y/(2n_-) ) \label{original-numer}\\ 
& = n \sin( (1 + 1/n) y) - (n+1) \sin( y). \label{simplified-numer}
\end{align}
A Taylor expansion of \eqref{simplified-numer} reveals that it is an alternating series with monotone coefficients.  Thus we can estimate error by the next term in the truncated series:
\begin{align*}
|n \sin( (1 + 1/n) y) - (n+1) \sin( y)| \leq \frac{|y|^3}{3!} (1+n)\left[ (1+1/n)^2 - 1 \right] \leq C |y^3|.
\end{align*}
Returning to \eqref{original-numer}, it is easy to see that because $|\sin (y/(2n_-))| \leq |y|/(2n_-)$, \eqref{original-numer} is bounded above by $C(1+|y|)$.  Combining these estimates, with \eqref{denom-bound} we have \eqref{Dprime-est}.  The preceding calculations applied to this situation show
\begin{align}
\|D'_n\|_{L^p(\epsilon/n_-,2\pi -\epsilon/n_-)} &\leq C_p(1+\epsilon) \left(\frac{n_-}{1+\epsilon}\right)^{2-1/p}.\label{derivative-estimates}
\end{align}
Again, setting $\epsilon = 0$ gives a special case of Theorem~\ref{Theorem:Galeev}.

\bibliographystyle{plain}
\bibliography{/Users/trogdon/Dropbox/References/library.bib}

\end{document}